\documentclass[preprint,12pt,1p]{elsarticle}
\makeatletter
\def\ps@pprintTitle{%
  \let\@oddhead\@empty
  \let\@evenhead\@empty
  \let\@oddfoot\@empty
  \let\@evenfoot\@oddfoot
}
\makeatother
\usepackage{graphics}
\usepackage{epsfig}
\usepackage{mathptmx}
\usepackage{amsmath}
\usepackage{amssymb}
\usepackage{hyperref}
\usepackage{fullpage}
\usepackage{amsthm}
\usepackage{lineno}

\theoremstyle{plain}
\newtheorem{theorem}{Theorem}[section]

\theoremstyle{definition}
\newtheorem{definition}{Definition}[section]
\theoremstyle{example}

\theoremstyle{remark}
\newtheorem{remark}{Remark}[section]
\newlength{\defbaselineskip}
\journal{Computer aided geometric design}
\begin{document}
\begin{frontmatter}
\title{TRANSVERSAL INTERSECTION CURVES OF HYPERSURFACES IN $\mathbb{R}^{5}$}
\author[nd]{Mohamd Saleem Lone}
\ead{saleemraja2008@gmail.com}
\author[nc]{O. Al\'{e}ssio}
\ead{osmar@mathematica.uftm.edu.br}
\author[nb]{Mohammad Jamali}
\ead{jamali\_dbd@yahoo.co.in}
\author[na]{Mohammad Hasan Shahid\corref{cor1}}
\ead{hasan\_jmi@yahoo.com}
\address[nd]{Central University of Jammu, Jammu, 180011, India.}
\address[nc]{Universidade Federal do Tri$\dot{a}$ngulo Mineiro-UFTM, Uberaba, MG, Brasil}
\address[nb]{Department of Mathematics, Al-Falah University, Haryana, 121004, India}
\address[na]{Department of Mathematics, Jamia Millia Islamia, New Delhi-110 025, India}
\cortext[cor1]{Corresponding author}
\begin{abstract}
In this paper we present the algorithms for calculating the differential geometric properties $\{t,n,b_{1},b_{2},b_{3},\kappa_{1},\kappa_{2},\kappa_{3},\kappa_{4}\}$, geodesic curvature and geodesic torsion of the transversal intersection curve of four hypersurfaces (given by parametric representation) in Euclidean space $\mathbb{R}^{5}$. In transversal intersection the normals of the surfaces at the intersection point are linearly independent, while as in nontransversal intersection the normals of the surfaces at the intersection point are linearly dependent.
\end{abstract}
\begin{keyword}
\texttt Hypersurfaces, transversal intersection, non-transversal intersection.
\end{keyword}
\end{frontmatter}

\section{Introduction}
The surface-surface intersection problem is a fundamental process needed in modeling shapes in CAD/CAM system. It is useful in the representation of the design of complex objects and animations.
The two types of surfaces most used in geometric designing are parametric and implicit surfaces. For that reason different methods have been given for either parametric-parametric or implicit-implicit surface intersection curves in $\mathbb{R}^{3}$. The numerical marching method is the most widely used method for computing the intersection  curves in $\mathbb{R}^{3}$ and $\mathbb{R}^{4}$. The marching method involves generation of sequences of points of an intersection curve in the direction prescribed by the local geometry(Bajaj et al., 1988; Patrikalakis, 1993). To compute the intersection curve with precision and efficiency, approaches of superior order are necessary, that is, they are needed to obtain the geometric properties of the intersection curves. While differential geometry of a parametric curve in $\mathbb{R}^{3}$ can be found in textbooks such as Struik(1950), Willmore (1959), Stoker (1969), Spivak (1975), do Carmo (1976), differential geometry of a parametric curve in $\mathbb{R}^{n}$ can be found in the textbook such as in klingenberg (1978) and in the contemporary literature on Geometric Modelling (Farin, 2002; Hoschek and Lasser 1993), but there is only a scarce of literature on the differential geometry of intersection curves. Willmore (1959) and Al\'{e}ssio (2006) presented algorithms to obtain the unit tangent, unit principle normal, unit binormal, curvature and torsion of the transversal intersection curve of two implicit surfaces. Hartmann (1996) provided formulas for computing the curvature of the intersection curves for all types of intersection problems in $\mathbb{R}^{3}$. Ye and Maekawa (1999) presented algorithms for computing the differential geometric properties of both transversal and tangential intersection curves of two surfaces. Al\'{e}ssio (2009) formulated the algorithms for obtaining the geometric properties of intersection curves of three implicit hypersurfaces in $\mathbb{R}^{4}$. Based on the work of Al\'{e}ssio (2009), Mustufa D\"{u}ld\"{u}l (2010) worked with three parametric hypersurfaces in  $\mathbb{R}^{4}$ to derive the algorithms for differential geometric properties of transversal intersection. Abdel-All et al. (2012) provided algorithms for geometric properties of implicit-implict-parametric and implicit-parametric-parametric hypersurfaces in $\mathbb{R}^{4}$. Naeim-Badr et al. (2014) obtained algorithms for differential geometric properties of non-transversal intersection curves of three parametric hypersurfaces in $\mathbb{R}^{4}$. Recently Naeim Badr, Abdel-All et al. (2015) derived the algorithms for non-transversal intersection curves of implicit-parametric-parametric and implicit-implicit-parametric hypersurfacres in $\mathbb{R}^{4}$. To obtain the first geodesic curvature ($\kappa_{1g}^{S_{i}}$) and the first geodesic torsion ($\tau_{1g}^{S_{i}}$) for the transversal intersection curve of $4$ parametric
hypersurfaces in $\mathbb{R}^{5}$, \ we need to derive the Darboux frame $\left\{  U_{1}^{M_{i}},\cdots,U_{5}^{M_{i}}\right\}  $. The Darboux frame is
obtained by using the Gram-Schmidt orthogonalization process.\newline
In this paper we extended the methods of Mustufa D\"{u}ld\"{u}l\cite{b18}, to obtain the Frenet frame $\{t,n,b_{1},b_{2},b_{3}\}$ and curvatures $\{\kappa_{1},\kappa_{2},\kappa_{3},\kappa_{4}\}$ of transversal intersection curve of four parametric hypersurfaces in $\mathbb{R}^{5}$. In section 2 we introduce some notations and reviews of the differential geometry of curves and surfaces in $\mathbb{R}^{5}$. In section 3 we find the formulas for computing the properties of transversal intersection of four parametric hypersurfaces in $\mathbb{R}^{5}$. In section 4 we derive the formulas for obtaining the geodesic curvature and geodesic torsion of the intersecting curve with respect to four hypersurfaces. Finally, to be more constructive we present an example in section 5. Moreover in addition to the use of classical results of differential geometry we will also make use of Matlab/Mathematica.

\section{ Preliminaries }
\begin{definition}Let $ e_{1},e_{2},e_{3},e_{4},e_{5}$ be the standard basis of five dimensional Euclidean space $E^{5}$. The vector product of the vectors $x=\sum_{i=1}^{5}x_{i}e_{i}$, $y=\sum_{i=1}^{5}y_{i}e_{i}$, $ z=\sum_{i=1}^{5}z_{i}e_{i}$ and $w=\sum_{i=1}^{5}w_{i}e_{i}$  is defined by
\begin{equation}\label{p1}x\otimes y\otimes z\otimes w =\left| \begin{array}{ccccc}
                                                            e_{1} & e_{2} & e_{3} & e_{4} & e_{5} \\
                                                            x_{1} & x_{2} & x_{3} & x_{4} & x_{5} \\
                                                            y_{1} & y_{2} & y_{3} & y_{4} & y_{5} \\
                                                            z_{1} & z_{2} & z_{3} & z_{4} & z_{5} \\
                                                            w_{1} & w_{2} & w_{3} & w_{4} & w_{5}
                                                          \end{array}
 \right|\end{equation} The vector product $x\otimes y \otimes z \otimes w$ yields a vector that is orthogonal to $x$, $y$, $z$, $w$.\end{definition}
 let $R\subset E^{5}$ be a regular hypersurface given by $\Phi = \Phi(u_{1},u_{2},u_{3},u_{4})$ and $\alpha : I\subset\mathbb{R}\rightarrow\Phi$ ba an arbitrary curve with arc length parametrisation. If $t,n,b_{1},b_{2},b_{3}$ is the Frenet Frame along $\alpha$
 \begin{eqnarray}\label{p2}\begin{split}t^{'}&=\kappa_{1}n&\\
 n^{'}&=-\kappa_{1}t+\kappa_{2}b_{1}&\\
 b_{1}^{'}&=-\kappa_{2}n+\kappa_{3}b_{2}&\\
 b_{2}^{'}&=-\kappa_{3}b_{1}+\kappa_{4}b_{3}&\\
 b_{3}^{'}&=-\kappa_{4}b_{2}&\end{split}\end{eqnarray}
 Where $t$, $n$, $b_{1}$, $b_{2}$ and $b_{3}$ denote the tangent, the principal normal, the first binormal, the second binormal and third binormal vector fields. The normal vector  $n$ is the normalised acceleration vector $\alpha^{''}$. The unit vector $b_{1}$ is determined such that $n^{'}$ can be decomposed into two components, a tangent one in the direction of $t$ and a normal one in the direction of $b_{1}$. The unit vector $b_{2}$ is determined such that $b_{1}^{'}$ can be decomposed into two components, a normal one and $b_{2}$. The unit vector $b_{3}$ is the unique unit vector field perpendicular to four dimensional subspace $\{t,n,b_{1},b_{2}\}$. The functions $\kappa_{1}$, $\kappa_{2}$, $\kappa_{3}$ and $\kappa_{4}$ are the first, second, third and fourth curvatures of $\alpha(s)$. The first, second, third and fourth curvatures measures how rapidly the curve pulls away in a neighbourhood of $s$, from the tangent line, from planar curve, from three dimensional curve and from the four dimensional curve at $s$, respectively.\newline
 Now, using the Frenet Frame we have the derivatives of $\alpha$ as
 \begin{eqnarray}\label{p3} \alpha^{'}&=&t, \hspace{1cm} \alpha^{''}=t^{'}=\kappa_{1}n,\hspace{1cm}\alpha^{'''}=-\kappa_{1}^{2}t+\kappa_{1}^{'}n+\kappa_{1}\kappa_{2}b_{1}\\
 \label{p4}\alpha^{(4)}&=&-3\kappa_{1}\kappa_{1}^{'}t + (-\kappa_{1}^{3}+\kappa_{1}^{''}-\kappa_{1}\kappa_{2}^{2})n+(2\kappa_{1}^{'}\kappa_{2}+\kappa_{1}\kappa_{2}^{'})b_{1}+\kappa_{1}\kappa_{2}\kappa_{3}b_{2}\\
 \label{p5}\alpha^{(5)}&=&(-3(\kappa_{1}^{'})^{2}-4\kappa_{1}k_{1}^{''}+\kappa_{1}^{4}+
 \kappa_{1}^{2}\kappa_{2}^{2})t+(-6\kappa_{1}^{2}\kappa_{1}^{'}+\kappa_{1}^{'''}-\kappa_{1}^{'}\kappa_{2}^{2}\nonumber\\
&& -3\kappa_{1}\kappa_{2}\kappa_{2}^{'}-2\kappa_{1}^{'}\kappa_{22}^{2})n+ (\kappa_{1}^{3}\kappa_{2}-\kappa_{1}\kappa_{2}^{3}+3\kappa_{1}^{''}\kappa_{2}+3\kappa_{1}^{'}\kappa_{2}^{'}+
\kappa_{1}\kappa_{2}^{''}-\kappa_{1}\kappa_{2}\kappa_{3}^{2})b_{1}\nonumber\\ &&(3\kappa_{1}^{'}\kappa_{2}\kappa_{3}+2\kappa_{1}\kappa_{2}^{'}\kappa_{3}+2\kappa_{1}\kappa_{2}^{'}\kappa_{3}+\kappa_{1}\kappa_{2}\kappa_{3}^{'})b_{2}+
\kappa_{1}\kappa_{2}\kappa_{3}\kappa_{4}b_{3}
\end{eqnarray}
Also since $\Phi$ is regular, the partial derivatives $\Phi_{1}$, $\Phi_{2}$, $\Phi_{3}$, $\Phi_{4}$, where $(\Phi_{i}=\frac{\partial \Phi}{\partial u_{i}})$ are linearly independent at every point of $\Phi$ , i.e., $\Phi_{1}\otimes \Phi_{2} \otimes \Phi_{3} \otimes \Phi_{4}\neq 0$. Thus, the unit normal vector of $\Phi$ is given by
$$N=\frac{\Phi_{1}\otimes \Phi_{2} \otimes \Phi_{3} \otimes \Phi_{4}}{\|\Phi_{1}\otimes \Phi_{2} \otimes \Phi_{3} \otimes \Phi_{4}\|}.$$
Furthermore, the first, second, and the third binormal vectors of the curve are given by
\begin{equation}\label{p6}b_{3}=\frac{\alpha^{'}\otimes \alpha^{''}\otimes \alpha^{'''}\otimes \alpha^{(4)}}{\|\alpha^{'}\otimes \alpha^{''}\otimes \alpha^{'''}\otimes \alpha^{(4)}\|}, \hspace{.5cm}b_{2}=\frac{b_{3} \otimes \alpha^{'}\otimes \alpha^{''} \otimes \alpha^{'''}}{\|b_{3} \otimes \alpha^{'}\otimes \alpha^{''} \otimes \alpha^{'''}\|}, \hspace{.5cm} b_{1} =\frac{ b_{2} \otimes b_{3}\otimes \alpha^{'} \otimes \alpha^{''}}{\| b_{2} \otimes b_{3}\otimes \alpha^{'} \otimes \alpha^{''}\|}
 \end{equation}
 and the curvatures are obtained with
 \begin{equation}\label{p7}\kappa_{1}=\|\alpha^{''}\|, \hspace{.3cm} \kappa_{2}=\frac{\langle\alpha^{'''},b_{1}\rangle}{\kappa_{1}}
 ,\hspace{.3cm} \kappa_{3}=\frac{\langle\alpha^{(4)},b_{2}\rangle}{\kappa_{1}\kappa_{2}}
 ,\hspace{.3cm} \kappa_{4}=\frac{\langle\alpha^{(5)},b_{3}\rangle}{\kappa_{1}\kappa_{2}\kappa_{3}}\end{equation}
 On the other hand, since the curve $\alpha(s)$ lies on $\Phi$, we may write\newline
 $\alpha(s)=\Phi(u_{1}(s),u_{2}(s),u_{3}(s),u_{4}(s)).$\newline
 We then have
 \begin{equation}\label{p8}\alpha^{'}(s)= \sum_{i=1}^{4}\Phi_{i} u_{i}^{'}\end{equation}
 \begin{equation}\label{p9}\alpha^{''}(s)= \sum_{i=1}^{4}\Phi_{i} u_{i}^{''}+\sum_{i,j=1}^{4}\Phi_{ij} u_{i}^{'}u_{j}^{'}\end{equation}
 \begin{equation}\label{p10}\alpha^{'''}(s)= \sum_{i=1}^{4}\Phi_{i} u_{i}^{'''}+3\sum_{i,j=1}^{4}\Phi_{ij} u_{i}^{''}u_{j}^{'}+\sum_{i,j,k=1}^{4}\Phi_{ijk} u_{i}^{'}u_{j}^{'}u_{k}^{'}\end{equation}
 \begin{eqnarray}\label{p11}\nonumber\alpha^{(4)}(s)=&& \sum_{i=1}^{4}\Phi_{i} u_{i}^{(4)}+4\sum_{i,j=1}^{4}\Phi_{ij} u_{i}^{'''}u_{j}^{'}+3\sum_{i,j=1}^{4}\Phi_{ij} u_{i}^{''}u_{j}^{''}+6\sum_{i,j,k=1}^{4}\Phi_{ijk} u_{i}^{''}u_{j}^{'}u_{k}^{'}+\\
 &&\sum_{i,j,k,l=1}^{4}\Phi_{ijkl} u_{i}^{'}u_{j}^{'}u_{k}^{'}u_{l}^{'}\end{eqnarray}
 \begin{eqnarray}\label{p12}\nonumber\alpha^{(5)}&=& \sum_{i=1}^{4}\Phi_{i} u_{i}^{(5)} + 5\sum_{i,j=1}^{4}\Phi_{ij}u_{i}^{(4)}u_{j}^{'}+10\sum_{i,j=1}^{4}\Phi_{ij}u_{i}^{'''}u_{j}^{''}+10\sum_{i,j,k=1}^{4}\Phi_{ijk}u_{i}^{'''}u_{j}^{'}u_{k}^{'}\\
 \nonumber&&+15\sum_{i,j,k=1}^{4}\Phi_{ijk}u_{i}^{''}u_{j}^{''}u_{k}^{'}+10\sum_{i,j,k,l=1}^{4}\Phi_{ijkj}u_{i}^{''}u_{j}^{'}u_{k}^{'}u_{l}^{'}\\
 &&+\sum_{i,j,k,l,r=1}^{4}\Phi_{ijklr}u_{i}^{'}u_{j}^{'}u_{k}^{'}u_{l}^{'}u_{r}^{'}\end{eqnarray}
 \section{The curvature of the transversal intersection of hypersurfaces}
 Let $R_{1}$, $R_{2}$, $R_{3}$ and $R_{4}$ be four regular transversally intersecting hypersurfaces given by parametric equations  $\Phi^{i}=\Phi^{i}(u_{1}^{i},u_{2}^{i},u_{3}^{i},u_{4}^{i})$, $(i=1,2,3,4)$. Then the unit normal vector of these hypersurfaces are
\begin{equation*}N_{i}=\frac{\Phi_{1}^{i} \otimes\Phi_{2}^{i}\otimes\Phi_{3}^{i}\otimes\Phi_{4}^{i}}{\|\Phi_{1}^{i} \otimes\Phi_{2}^{i}\otimes\Phi_{3}^{i}\otimes\Phi_{4}^{i}\|},\hspace{1cm} i=1,2,3,4\end{equation*}
Since the intersection is transversal, the normal vectors of these hypersurfaces at the intersection points are linearly independent, i.e., $N_{1}\otimes N_{2}\otimes N_{3}\otimes N_{4}\neq 0$. It is assumed that the intersection is a smooth curve say $\alpha$, in $E^{5}$. Let the intersection curve $\alpha$ be parameterised by arc lenght function $s$. Then, at the intersection point $\alpha(s_{0})=P$, the unit tangent vector $t$ of the intersection curve $\alpha$ can be found by the vector product of the normal vectors at $P$.
\begin{equation}\label{p13}t=\frac{N_{1}\otimes N_{2} \otimes N_{3} \otimes N_{4}}{\|N_{1}\otimes N_{2} \otimes N_{3} \otimes N_{4}\|}\end{equation}
\subsection{First curvature of the transversal intersection curve}
Now, we find the first curvature of the intersection curve at $P$. Since $t^{'}$ is orthogonal to $t$, we may write
\begin{equation}\label{p14}\alpha^{''}=t^{'}=\sum_{i=1}^{4}a_{i}N_{i}, \hspace{1cm} a_{i}\in \mathbb{R}, i=1,2,3,4\end{equation}
Thus, we need to calculate the scalars $a_{i}$ to find $\alpha^{''}$.
If we take the dot product of both hand sides of $(\ref{p14})$ with $N_{i}$, we have
\begin{eqnarray}\label{p15}\begin{split}
\langle N_{1}, N_{1}\rangle a_{1}+\langle N_{2}, N_{1}\rangle a_{2}+\langle N_{3}, N_{1}\rangle a_{3}+\langle N_{4}, N_{1}\rangle a_{4}=\kappa_{n}^{1},\\
\langle N_{1}, N_{2}\rangle a_{1}+\langle N_{2}, N_{2}\rangle a_{2}+\langle N_{3}, N_{2}\rangle a_{3}+\langle N_{4}, N_{2}\rangle a_{4}=\kappa_{n}^{2},\\
\langle N_{1}, N_{3}\rangle a_{1}+\langle N_{2}, N_{3}\rangle a_{2}+\langle N_{3}, N_{3}\rangle a_{3}+\langle N_{4}, N_{3}\rangle a_{4}=\kappa_{n}^{3},\\
\langle N_{1}, N_{4}\rangle a_{1}+\langle N_{2}, N_{4}\rangle a_{2}+\langle N_{3}, N_{4}\rangle a_{3}+\langle N_{4}, N_{4}\rangle a_{4}=\kappa_{n}^{4}\end{split}\end{eqnarray}
Where, $\kappa_{n}^{i}=\langle t^{'}, N_{i}\rangle, \hspace{.5cm}(i=1,2,3,4)$ and $\langle, \rangle$ is the scalar product.
Hence we must compute $\kappa_{n}^{1}$, $\kappa_{n}^{2}$, $\kappa_{n}^{3}$, $\kappa_{n}^{4}$ at $P$ to find the scalars $a_{i}$
On using $(\ref{p9})$, we obtain
\begin{eqnarray}\label{p17}\begin{split}\kappa_{n}^{i}=&\Pi_{11}^{i}({u_{1}^{i}}^{'})^{2}+\Pi_{22}^{i}({u_{2}^{i}}^{'})^{2}+
\Pi_{33}^{i}({u_{3}^{i}}^{'})^{2}+\Pi_{44}^{i}({u_{4}^{i}}^{'})^{2}+2(\Pi_{12}^{i}{u_{1}^{i}}^{'}{u_{2}^{i}}^{'}+
\Pi_{13}^{i}{u_{1}^{i}}^{'}{u_{3}^{i}}^{'}&
\\&+\Pi_{14}^{i}{u_{1}^{i}}^{'}{u_{4}^{i}}^{'}+\Pi_{23}^{i}{u_{2}^{i}}^{'}{u_{3}^{i}}^{'}
+\Pi_{24}^{i}{u_{2}^{i}}^{'}{u_{4}^{i}}^{'}+\Pi_{34}^{i}{u_{3}^{i}}^{'}{u_{4}^{i}}^{'})& \end{split}\end{eqnarray}
Where $\Pi_{lm}^{i}$, $l,m = 1,2,3,4$ are the second fundamental form coefficients of the hypersurfaces $\Phi^{i}$.
Since the unit tangent is known from $(\ref{p13})$ and $\|\Phi_{1}^{i}\otimes \Phi_{2}^{i}\otimes \Phi_{3}^{i}\otimes \Phi_{4}^{i}\| \neq 0$, the scalar multiplication of both hand sides of $(\ref{p8})$ with $\Phi_{1}^{i}$, $\Phi_{2}^{i}$, $\Phi_{3}^{i}$ and $\Phi_{4}^{i}$, respectively yields a linear system of four equations as
\begin{eqnarray}\label{p18}\begin{split}
\Upsilon_{11}^{i}{u_{1}^{i}}^{'}+\Upsilon_{12}^{i}{u_{2}^{i}}^{'}+\Upsilon_{13}^{i}{u_{3}^{i}}^{'}+\Upsilon_{14}^{i}{u_{4}^{i}}^{'}=\langle t,\Phi_{1}^{i}\rangle\\
\Upsilon_{21}^{i}{u_{1}^{i}}^{'}+\Upsilon_{22}^{i}{u_{2}^{i}}^{'}+\Upsilon_{23}^{i}{u_{3}^{i}}^{'}+\Upsilon_{24}^{i}{u_{4}^{i}}^{'}=\langle t,\Phi_{2}^{i}\rangle\\
\Upsilon_{31}^{i}{u_{1}^{i}}^{'}+\Upsilon_{32}^{i}{u_{2}^{i}}^{'}+\Upsilon_{33}^{i}{u_{3}^{i}}^{'}+\Upsilon_{34}^{i}{u_{4}^{i}}^{'}=\langle t,\Phi_{3}^{i}\rangle\\
\Upsilon_{41}^{i}{u_{1}^{i}}^{'}+\Upsilon_{42}^{i}{u_{2}^{i}}^{'}+\Upsilon_{43}^{i}{u_{3}^{i}}^{'}+\Upsilon_{44}^{i}{u_{4}^{i}}^{'}=\langle t,\Phi_{4}^{i}\rangle
\end{split}\end{eqnarray}
with respect to $u_{1}^{'}$, $u_{2}^{'}$, $u_{3}^{'}$ and $u_{4}^{'}$ where $\Upsilon_{lm}^{i}, \hspace{.2cm}l,m=1,2,3,4$ are the first fundamental form coefficients of the hypersurface $\Phi^{i}$. Substituting the solutions of these systems into $(\ref{p17})$ gives us $\kappa_{n}^{i}, \hspace{.5cm}i=1,2,3,4$. Then using matlab/mathematica the system of linear equations in (\ref{p15}) can be solved for $a_{1}$, $a_{2}$, $a_{3}$ and $a_{4}$.Thus the first curvature of the intersection curve at $P$ is obtained from (\ref{p14}) and the first equation of (\ref{p7}).
\remark{If the normal vectors are mutually orthogonal to each other at the intersection point, then the first curvature is given by $\kappa_{1}=\sqrt{a_{1}^{2}+a_{2}^{2}+a_{3}^{2}+a_{4}^{2}}$}.
\subsection{Second curvature}
To compute the second curvature we have to find the third derivative of the intersection curve $\alpha$.\newline
Since $N_{i}, i=1,2,3,4$ are orthogonal to $t$, we may write $c_{1}N_{1}+c_{2}N_{2}+c_{3}N_{3}+c_{4}N_{4}$ instead of terms $\kappa_{1}^{'}n + \kappa_{1}\kappa_{2}b_{1}$ in $\alpha^{'''}$, i.e.,
\begin{equation}\label{p20}\alpha^{'''}=-\kappa_{1}^{2}t+c_{1}N_{1}+c_{2}N_{2}+c_{3}N_{3}+c_{4}N_{4}.\end{equation}
Taking the dot product of both hand sides of above equation with $N_{i}$, we obtain
\begin{eqnarray}\label{p21}\begin{split}
\langle N_{1}, N_{1}\rangle c_{1}+\langle N_{2}, N_{1}\rangle c_{2}+\langle N_{3}, N_{1}\rangle c_{3}+\langle N_{4}, N_{1}\rangle c_{4}=\mu_{1}\\
\langle N_{1}, N_{2}\rangle c_{1}+\langle N_{2}, N_{2}\rangle c_{2}+\langle N_{3}, N_{2}\rangle c_{3}+\langle N_{4}, N_{2}\rangle c_{4}=\mu_{2}\\
\langle N_{1}, N_{3}\rangle c_{1}+\langle N_{2}, N_{3}\rangle c_{2}+\langle N_{3}, N_{3}\rangle c_{3}+\langle N_{4}, N_{3}\rangle c_{4}=\mu_{3}\\
\langle N_{1}, N_{4}\rangle c_{1}+\langle N_{2}, N_{4}\rangle c_{2}+\langle N_{3}, N_{4}\rangle c_{3}+\langle N_{4}, N_{4}\rangle c_{4}=\mu_{4}\end{split}\end{eqnarray}
Where $\mu_{i}=\langle\alpha^{'''}, N_{i} \rangle,i=1,2,3,4$ and $c_{i}\in \mathbb{R}$.
Now, let us find the unknown scalars $\mu_{i}$. Using $(\ref{p10})$, we have
\begin{eqnarray}\label{p23}\begin{split}
\mu_{r}= 3\sum_{i,j=1}^{4}\langle\Phi_{ij}^{r}, N_{r}\rangle u_{i}^{''}u_{j}^{'}+\sum_{i,j,k=1}^{4}\langle\Phi_{ijk}^{r}, N_{r}\rangle u_{i}^{'}u_{j}^{'}u_{k}^{'}
\end{split}\end{eqnarray}
Since the components $u_{i}^{'}$ are known from the system $(\ref{p18})$, we have to find $u_{i}^{''}$. To obtain $u_{i}^{''}$
we use $(\ref{p9})$ and write
\begin{equation}\label{p24}\Delta_{r}=\sum_{i,j=1}^{4}\Phi_{ij}^{r}u_{i}^{'}u_{j}^{'}, \hspace{1cm}r=1,2,3,4\end{equation}Then we have
\begin{eqnarray}\label{p25}\begin{split}
\Upsilon_{11}^{i}{u_{1}^{i}}^{''}+\Upsilon_{12}^{i}{u_{2}^{i}}^{''}+\Upsilon_{13}^{i}{u_{3}^{i}}^{''}+\Upsilon_{14}^{i}{u_{4}^{i}}^{''}=\langle t^{'}-\Delta_{i},\Phi_{1}^{i}\rangle\\
\Upsilon_{21}^{i}{u_{1}^{i}}^{''}+\Upsilon_{22}^{i}{u_{2}^{i}}^{''}+\Upsilon_{23}^{i}{u_{3}^{i}}^{''}+\Upsilon_{24}^{i}{u_{4}^{i}}^{''}=\langle t^{'}-\Delta_{i},\Phi_{2}^{i}\rangle\\
\Upsilon_{31}^{i}{u_{1}^{i}}^{''}+\Upsilon_{32}^{i}{u_{2}^{i}}^{''}+\Upsilon_{33}^{i}{u_{3}^{i}}^{''}+\Upsilon_{34}^{i}{u_{4}^{i}}^{''}=\langle t^{'}-\Delta_{i},\Phi_{3}^{i}\rangle\\
\Upsilon_{41}^{i}{u_{1}^{i}}^{''}+\Upsilon_{42}^{i}{u_{2}^{i}}^{''}+\Upsilon_{43}^{i}{u_{3}^{i}}^{''}+\Upsilon_{44}^{i}{u_{4}^{i}}^{''}=\langle t^{'}-\Delta_{i},\Phi_{4}^{i}\rangle
\end{split}\end{eqnarray}
Which gives us required derivatives. Thus from (\ref{p23}), we find the values of $\mu_{i}$, which finally helps us to find the value of $c_{i}$ in system (\ref{p21}).
Thus the second curvature can be found from the second equation of (\ref{p7}), untill we find $b_{1}$.
\newline On using $(\ref{p7})$ we obtain
$\kappa_{1}^{'}=\langle \alpha^{'''}, n\rangle$.
\subsection{Third curvature}
To find the third curvature, we need to find the fourth derivative of the intersection curve $\alpha$ at $P$. \newline
Since $N_{i}$ is orthogonal to $t$, we may write $d_{1}N_{1}+d_{2}N_{2}+d_{3}N_{3}+d_{4}N_{4}$ instead of$(-\kappa_{1}^{3}+\kappa_{1}^{''}-\kappa_{1}\kappa_{2}^{2})n+(2\kappa_{1}^{'}\kappa_{2}+\kappa_{1}\kappa_{2}^{'})b_{1}+\kappa_{1}\kappa_{2}
\kappa_{3}b_{2}$ in $\alpha^{(4)}$, i.e.,
\begin{equation}\label{p26}\alpha^{(4)}=-3\kappa_{1}\kappa_{1}^{'}t+d_{1}N_{1}+d_{2}N_{2}+d_{3}N_{3}+d_{4}N_{4}\end{equation}
Taking the dot product of (\ref{p26}) with $N_{1}$, $N_{2}$, $N_{3}$ and $N_{4}$, we get
\begin{eqnarray}\begin{split}
\langle N_{1}, N_{1}\rangle d_{1}+\langle N_{2}, N_{1}\rangle d_{2}+\langle N_{3}, N_{1}\rangle d_{3}+\langle N_{4}, N_{1}\rangle d_{4}=\xi_{1}\\
\langle N_{1}, N_{2}\rangle d_{1}+\langle N_{2}, N_{2}\rangle d_{2}+\langle N_{3}, N_{2}\rangle d_{3}+\langle N_{4}, N_{2}\rangle d_{4}=\xi_{2}\\
\langle N_{1}, N_{3}\rangle d_{1}+\langle N_{2}, N_{3}\rangle d_{2}+\langle N_{3}, N_{3}\rangle d_{3}+\langle N_{4}, N_{3}\rangle d_{4}=\xi_{3}\\
\langle N_{1}, N_{4}\rangle d_{1}+\langle N_{2}, N_{4}\rangle d_{2}+\langle N_{3}, N_{4}\rangle d_{3}+\langle N_{4}, N_{4}\rangle d_{4}=\xi_{4}\end{split}\end{eqnarray}
Where $\xi_{i}=\langle\alpha^{(4)}, N_{i} \rangle,i=1,2,3,4$ and $d_{i}\in \mathbb{R}$\newline
Now to find $d_{i}$ we have to find $\xi_{i}$. For that taking the dot product of $\alpha^{(4)}$ with $N_{i}$, we obtain
\begin{eqnarray}\label{p29}\nonumber\xi_{r}=&& 4\sum_{i,j=1}^{4}\langle\Phi_{ij}^{r}, N_{r}\rangle u_{i}^{'''}u_{j}^{'}+3\sum_{i,j=1}^{4}\langle\Phi_{ij}^{r},N_{r}\rangle u_{i}^{''}u_{j}^{''}+6\sum_{i,j,k=1}^{4}\langle \Phi_{ijk}^{r},N_{r}\rangle u_{i}^{''}u_{j}^{'}u_{k}^{'}+\\
 &&\sum_{i,j,k,l=1}^{4}\langle\Phi_{ijkl}^{r},N_{r}\rangle u_{i}^{'}u_{j}^{'}u_{k}^{'}u_{l}^{'}\end{eqnarray}
 Since $u_{i}^{'}$, $u_{i}^{''}$ are already known, so to find $\xi_{r}$ we have to find $u_{i}^{'''}$. These derivatives can be found by taking the product of both hand side of (\ref{p10}) with $\Phi_{1}^{i}$, $\Phi_{2}^{i}$, $\Phi_{3}^{i}$ and $\Phi_{4}^{i}$, respectively.Hence, we can compute the Frenet vectors at $P$ of the intersection curve by finding $b_{3}$, $b_{2}$ and $b_{1}$ - the third, second and first binormal vectors from the equations in (\ref{p6}) as now $\alpha^{'}$, $\alpha^{''}$, $\alpha^{'''}$ and $\alpha^{(4)}$ are at our disposal.
Thus on using the binormal vector $b_{1}$ and $(\ref{p20})$, the second curvature of the intersection curve at $P$ is obtained from the second equation of $(\ref{p7})$.
 Since $\kappa_{1}$, $\kappa_{2}$ and $b_{2}$ are already known, the third curvature can be now found from the third equation of (\ref{p7}).
 \subsection{Fourth curvature}
To obtain the forth curvature $\kappa_{4}$, we need to find the fifth derivative of the intersection curve of $\alpha$ at $P$ .
Similar to third and fourth derivative of the curve $\alpha$, we may write
\begin{equation}\label{p30}\alpha^{(5)}=
\{-3({\kappa_{1}^{'}})^{2}-4\kappa_{1}\kappa_{2}^{''}+\kappa_{1}^{(4)}+\kappa_{1}^{2}\kappa_{2}^{2}\}t+m_{1}N_{1}+m_{2}N_{2}+m_{3}N_{4}+n_{4}N_{4}
\end{equation}
Where, the system of equations for unknowns is
\begin{eqnarray}\begin{split}
\langle N_{1}, N_{1}\rangle m_{1}+\langle N_{2}, N_{1}\rangle m_{2}+\langle N_{3}, N_{1}\rangle m_{3}+\langle N_{4}, N_{1}\rangle m_{4}=\eta_{1}\\
\langle N_{1}, N_{2}\rangle m_{1}+\langle N_{2}, N_{2}\rangle m_{2}+\langle N_{3}, N_{2}\rangle m_{3}+\langle N_{4}, N_{2}\rangle m_{4}=\eta_{2}\\
\langle N_{1}, N_{3}\rangle m_{1}+\langle N_{2}, N_{3}\rangle m_{2}+\langle N_{3}, N_{3}\rangle m_{3}+\langle N_{4}, N_{3}\rangle m_{4}=\eta_{3}\\
\langle N_{1}, N_{4}\rangle m_{1}+\langle N_{2}, N_{4}\rangle m_{2}+\langle N_{3}, N_{4}\rangle m_{3}+\langle N_{4}, N_{4}\rangle m_{4}=\eta_{4}\end{split}\end{eqnarray}
and $\eta_{i}=\langle\alpha^{(5)}, N_{i}\rangle$.
Projecting (\ref{p12}) onto the unit vector $N_{i}$, respectively, we obtain $\eta_{i}$ depending on $u_{i}^{(4)}$ besides $u_{i}^{'''}$, $u_{2}^{''}$ and $u_{1}^{'}$. Except $u_{i}^{(4)}$ all are known. So to find $u_{i}^{(4)}$ taking the scalar product of (\ref{p11}) with $\Phi_{1}^{i}$, $\Phi_{2}^{i}$, $\Phi_{3}^{i}$, $\Phi_{4}^{i}$, respectively.
Consequently, the fourth curvature of the intersection curve can be found from the last equation of (\ref{p7}).
\section{\textbf{Darboux Frame, First Geodesic Curvature and First, Second and Third Geodesic Torsion.}}
In this section, we derive the
Darboux frame $\left\{  U_{1}^{i},U_{2}^{i},U_{3}^{i},U_{4}^{i},U_{5}^{i}\right\}  ,$ the
first geodesic curvature ($\kappa_{1g}^{i}$) and the first, second and third geodesic
torsion ($\tau_{jg}^{i}$), $j=1,2,3$ for the transversal intersection curve of $4$ parametric hypersurfaces in $\mathbb{R}^{5}$.
\begin{definition}
Let $M_{i}$ be a regular hypersurface in $\mathbb{R}^{5}$ and $\alpha$ be a curve on $M_{i}$.
Then for each $i$, $1\leq i\leq 4,$ the function \newline$\kappa
_{ig}:I\longrightarrow\mathbb{R}$ \ \ \ \ \ defined for $s\in I$ by
\newline$\kappa_{ig}(s)=\left\langle U_{i}^{\prime}(s),U_{i+1}(s)\right\rangle
$ \newline\noindent is called the \textbf{i}$^{th}$\textbf{\ geodesic
curvature function} of the curve $\alpha$ and $\kappa_{ig}(s)$ is called the
\textbf{i}$^{th}$\textbf{\ geodesic curvature} of the curve $\alpha$ at
$\alpha(s).$
\end{definition}

\noindent

\noindent

\label{sec:computer}
\subsubsection{Darboux Frame:}
We are able to obtain a natural frame for the intersection curve-hypersurface
pair $\left(  \alpha(s),\mathbf{M}_{i}\right)  $, i.e., the frame $\left\{  U_{1}^{i},U_{2}^{i},U_{3}^{i},U_{4}^{i},U_{5}^{i}\right\}  $ by using the Gram-Schmidt orthogonalization process. By assumption the sets $\mathbf{\alpha}^{\prime}\left(  s\right)  ,$
$\mathbf{\alpha}^{\prime\prime}\left(  s\right)  ,$ $\mathbf{\alpha}^{\prime\prime\prime}\left(
s\right)  $ and $\mathbf{\alpha}^{(4)}\left(  s\right),$ and $\{N_{1}(p),N_{2}(p),N_{3}(p),N_{4}(p)\}$
are linearly independent.

 \noindent Fixing $U_{1}^{i}=\mathbf{\alpha}^{\prime}(p)$ and $U_{5}^{i}=N_{i},$ we have

\noindent The natural frame (Darboux frame) $\left\{  U_{1}^{i},U_{2}^{i},U_{3}^{i},U_{4}^{i},U_{5}^{i}\right\}  $ is
obtained, with $1\leq i\leq4.$

\bigskip\noindent{\small $
\begin{array}[c]{cc}
U_{1}^{i}= & \mathbf{\alpha}^{\prime}\left(  p\right) \\
U_{5}^{i}= & N_{i} \\
U_{2}^{i}= & \frac{-\left\langle \mathbf{\alpha}^{\prime\prime}(p),N_{i}\right\rangle N_{i}-\left\langle \mathbf{\alpha}^{\prime\prime}(p),U_{1}^{i}\right\rangle U_{1}^{i}+\mathbf{\alpha}^{\prime\prime}(p)}{\left\Vert -\left\langle \mathbf{\alpha}^{\prime\prime}(p),N_{i}\right\rangle N_{i}-\left\langle \mathbf{\alpha}^{\prime\prime}(p),U_{1}^{i}\right\rangle U_{1}^{i}+\mathbf{\alpha}^{\prime\prime}(p)\right\Vert },\\

U_{3}^{i}= & \frac{-\left\langle \mathbf{\alpha}^{\prime\prime\prime}(p),N_{i}\right\rangle N_{i}-\left\langle\mathbf{\alpha}^{\prime\prime\prime}(p),U_{2}^{i}\right\rangle U_{2}^{i}-\left\langle \mathbf{\alpha}^{\prime\prime\prime}(p),U_{1}^{i}\right\rangle U_{1}^{i}+\mathbf{\alpha}^{\prime\prime\prime}(p)}{\left\Vert -\left\langle \mathbf{\alpha}^{\prime\prime\prime}(p),N_{i}\right\rangle N_{i}-\left\langle\mathbf{\alpha}^{\prime\prime\prime}(p),U_{2}^{i}\right\rangle U_{2}^{i}-\left\langle \mathbf{\alpha}^{\prime\prime\prime}(p),U_{1}^{i}\right\rangle U_{1}^{i}+\mathbf{\alpha}^{\prime\prime\prime}(p\right\Vert },\\

U_{4}^{i}= & \frac{-\left\langle \mathbf{\alpha}^{(4)}(p),N_{i}\right\rangle N_{i}-\left\langle \mathbf{\alpha}^{(4)}(p),U_{3}^{i}\right\rangle U_{3}^{i}-\left\langle
\mathbf{\alpha}^{(4)}(p),U_{2}^{i}\right\rangle U_{2}^{i}-\left\langle \mathbf{\alpha}^{(4)}(p),U_{1}^{i}\right\rangle
U_{1}^{i}+\mathbf{\alpha}^{(4)}(p)}{\left\Vert -\left\langle \mathbf{\alpha}^{(4)}(p),N_{i}\right\rangle N_{i}-\left\langle \mathbf{\alpha}^{(4)}(p),U_{3}^{i}\right\rangle U_{3}^{i}-\left\langle
\mathbf{\alpha}^{(4)}(p),U_{2}^{i}\right\rangle U_{2}^{i}-\left\langle \mathbf{\alpha}^{(4)}(p),U_{1}^{i}\right\rangle
U_{1}^{i}+\mathbf{\alpha}^{(4)}(p)\right\Vert }.
\end{array}
$ }

\noindent The j-th geodesic curvature $\kappa_{jg}^{i}$ associated with i-th hypersurface is
\begin{equation}%
\begin{array}
[c]{ccc}%
\kappa_{jg}^{i} & = & \left\langle \left(  U_{j}^{i}\right)  ^{\prime
},U_{j+1}^{i}\right\rangle =-\left\langle \left(  U_{j+1}^{i}\right)
^{\prime},U_{j}^{i}\right\rangle, j \in \{1,2,3\} .
\end{array}
\label{eq:curvaturageodesica}%
\end{equation}

\noindent The j-th geodesic torsion $\tau_{jg}^{i}$ associated with i-th hypersurface  is
\begin{equation}%
\begin{array}
[c]{ccc}%
\tau_{jg}^{i} & = & -\left\langle \left(  U_{5}^{i}\right)  ^{\prime
},U_{j+1}^{i}\right\rangle =\left\langle \left(  U_{j+1}^{i}\right)
^{\prime},U_{5}^{i}\right\rangle, j \in \{1,2,3\}.
\end{array}
\label{eq:torsiongeodesica}%
\end{equation}
\subsubsection{\textbf{First geodesic and j-th torsion geodesic Formulas}}

\begin{theorem}
First Geodesic Curvature and the j-th geodesic torsion of the intersection curve of 4 parametric
hypersurfaces is
\begin{equation}
\kappa_{1g}^{i}=\sum_{j=1}^{4}a_j\left\langle N_{j},U_{2}^{i}\right\rangle
\label{eq:kg01}
\end{equation}
\begin{equation}
\tau_{jg}^{i}=-\left\langle \dfrac{(\overline{N}_{i})_{u_1}u^{\prime}_1 +(\overline{N}_{i})_{u_2}u^{\prime}_2+(\overline{N}_{i})_{u_3}u^{\prime}_3+(\overline{N}_{i})_{u_4}u^{\prime}_4}{\left\Vert \Phi^{i}_{1} \times \Phi^{i}_{2}\times \Phi^{i}_{3} \times \Phi^{i}_{4}\right\Vert } ,U_{j+1}^{i}\right\rangle , j=1,2,3.
\end{equation}

\end{theorem}

\noindent\textbf{Proof.}

For First Geodesic Curvature
\begin{equation}%
\begin{array}
[c]{ccc}%
\kappa_{1g}^{i} & = & \left\langle \left(  U_{1}^{i}\right)  ^{\prime
},U_{2}^{i}\right\rangle =\left\langle \mathbf{\alpha}^{\prime\prime}(s),U_{2}^{i}\right\rangle.
\end{array}
\end{equation}

Now using (14), Eq. (30) follows easily.

For the j-th geodesic torsion, we need derivative of $U_{5}^{i}=N_{i}=\dfrac{\Phi^{i}_{1} \times \Phi^{i}_{2}\times \Phi^{i}_{3} \times \Phi^{i}_{4}}{\left\Vert \Phi^{i}_{1} \times \Phi^{i}_{2}\times \Phi^{i}_{3} \times \Phi^{i}_{4}\right\Vert }$. If defined

$N_{i}=\dfrac{\overline{N}_{i}}{\left\|\overline{N}_{i}\right\|}, \ where \ \overline{N}_{i}=\Phi^{i}_{1}\times \Phi^{i}_{2}\times \Phi^{i}_{3} \times \Phi^{i}_{4} $. Hence we derive

$\frac{d}{ds}\left(U_{5}^{i}\left(s\right)  \right)=\dfrac{d N_{i}}{ds}=\dfrac{\dfrac{d \overline{N_{i}}}{ds}\left\|\overline{N}_{i}\right\|-\overline{N_{i}} \dfrac{d\left(\left\|\overline{N}_{i}\right\|\right)}{ds}}{{\left\|\overline{N}_{i}\right\|^2}}$

$\tau_{jg}^{i}=-\left\langle \dfrac{d N_{i}}{ds} ,U_{j+1}^{i}\right\rangle , j=1,2,3.$

$\tau_{jg}^{i}=-\left\langle \dfrac{\frac{d \overline{N_{i}}}{ds}\left\|\overline{N}_{i}\right\|}{{\left\|\overline{N}_{i}\right\|^2}}-\dfrac{\overline{N_{i}} \dfrac{d\left(\left\|\overline{N}_{i}\right\|\right)}{ds}}{{\left\|\overline{N}_{i}\right\|^2}} ,U_{j+1}^{i}\right\rangle , j=1,2,3.$

$\tau_{jg}^{i}=-\left\langle \dfrac{\frac{d \overline{N_{i}}}{ds}\left\|\overline{N}_{i}\right\|}{{\left\|\overline{N}_{i}\right\|^2}},U_{j+1}^{i}\right\rangle , j=1,2,3.$, because $ \left\langle \overline{N}_{i}, U_{j+1}^{i}\right\rangle =0, j=1,2,3.$

$\tau_{jg}^{i}=-\left\langle \dfrac{\frac{d \overline{N_{i}}}{ds}}{\left\|\overline{N}_{i}\right\|},U_{j+1}^{i}\right\rangle, j=1,2,3. $

$$
\frac{d \overline{N_{i}}}{ds}=(\overline{N}_{i})_{u_1}u_1^{\prime}+(\overline{N}_{i})_{u_2}u_2^{\prime}+ (\overline{N}_{i})_{u_3}u_3^{\prime}+ (\overline{N}_{i})_{u_4}{u_4}^{\prime}
$$

$
\begin{array}{cc}
(\overline{N}_{i})_{u_i}= \Phi^{i}_{1i} \times \Phi^{i}_{2}\times \Phi^{i}_{3} \times \Phi^{i}_{4}+\Phi^{i}_{1} \times \Phi^{i}_{2i}\times \Phi^{i}_{3} \times \Phi^{i}_{4}+ \\
\Phi^{i}_{1} \times \Phi^{i}_{2}\times \Phi^{i}_{3i} \times \Phi^{i}_{4}+\Phi^{i}_{1} \times \Phi^{i}_{2}\times \Phi^{i}_{3} \times \Phi^{i}_{4i}
\end{array}
.$
\noindent \section{Example}
Let $M_{1}$, $M_{2}$, $M_{3}$ and $M_{4}$ be the hypersurfaces given by, respectively
\begin{equation*}X^{1}(u_{1},u_{2},u_{3},u_{4})=(\sqrt{2}sinu_{3}\,cosu_{2}\,cosu_{1},\sqrt{2}sinu_{3}\,sinu_{2}\,sinu_{1},\sqrt{2}sinu_{3}\,cosu_{2},\sqrt{2}\,u_{4}\,cosu_{3},\frac{1}{\sqrt{2}}\,sinu_{3})\end{equation*}
\begin{equation*}X^{2}(u_{1},u_{2},u_{3},u_{4})=(u_{3}\,cosu_{1}\,cosu_{2},u_{3}\,sinu_{1}\,cosu_{2},u_{3}\,sinu_{2},u_{3},cosu_{2}\,cosu_{4})\end{equation*}
\begin{equation*}X^{3}(u_{1},u_{2},u_{3},u_{4})=(u_{3}\,cosu_{1}\,cosu_{2},u_{3}sinu_{1}\,cosu_{2},u_{4}\,sinu_{1},u_{3},sinu_{1}\,sinu_{2})\end{equation*}
\begin{equation*}X^{4}(u_{1},u_{2},u_{3},u_{4})=(\frac{1}{2}+\frac{1}{2}\,cosu_{1},\frac{1}{2}\,sinu_{1},u_{2},u_{3},\frac{u_{4}}{2})\end{equation*}
let us find the Frenet vectors and the curvatures of the intersection curve at the intersection point\newline
\begin{equation*}p=X^{1}\left(\frac{\pi}{4},\frac{\pi}{4},\frac{\pi}{4},1\right)=X^{2}\left(\frac{\pi}{4},\frac{\pi}{4},1,\frac{\pi}{4}\right)=
X^{3}\left(\frac{\pi}{4},\frac{\pi}{4},1,1\right)=X^{4}\left(\frac{\pi}{2},\frac{\sqrt{2}}{2},1,1\right)=\left(\frac{1}{2},
\frac{1}{2},\frac{\sqrt{2}}{2},1,\frac{1}{2}\right)\end{equation*}
The unit normals of these hypersurfaces are\newline
\begin{eqnarray*}N_{1}=\left(\frac{1}{\sqrt{6}},\frac{1}{\sqrt{6}},0,0,-\sqrt{\frac{2}{3}}\right),N_{2}=\left(-\frac{1}{2\sqrt{2}},-\frac{1}{2\sqrt{2}}
,-\frac{1}{2},\frac{1}{\sqrt{2}},0\right)\\
,N_{3}=\left(-\frac{2}{3},0,0,\frac{1}{3},-\frac{2}{3}\right),N_{4}=(0,1,0,0,0)\end{eqnarray*}
The non-vanishing first fundamental coefficients are\newline
\begin{eqnarray*}g_{11}^{1}=\frac{1}{2}\hspace{.5cm}g_{22}^{1}=1\hspace{.5cm}g_{33}^{1}=\frac{9}{4}\hspace{.5cm}g_{44}^{1}=1
\hspace{.5cm}g_{12}^{1}=\frac{1}{2}\hspace{.5cm}g_{23}^{1}=-\frac{1}{2}\hspace{.5cm}g_{34}^{1}=-1\\
g_{11}^{2}=\frac{1}{2}, \hspace{.5cm}g_{22}^{2}=\frac{5}{4}\hspace{.5cm}g_{33}^{2}=2,\hspace{.5cm}g_{44}^{2}=\frac{1}{4},\hspace{.5cm}g_{24}^{2}=\frac{1}{4},\hspace{.5cm}
g_{11}^{3}=\frac{5}{4},
\hspace{.5cm}g_{22}^{3}=\frac{3}{4},\\
g_{33}^{3}=\frac{3}{2},\hspace{.5cm}g_{44}^{3}=\frac{1}{2}, \hspace{.5cm}g_{12}^{3}=\frac{1}{4},\hspace{.5cm}g_{14}^{3}=\frac{1}{2},\hspace{.5cm}g_{23}^{3}=-\frac{1}{2},\hspace{.5cm}g_{11}^{4}=\frac{1}{4},
g_{22}^{4}=1,\\
g_{33}^{4}=1,\hspace{.5cm}g_{44}^{4}=\frac{1}{4}\end{eqnarray*}
The unit tangent at the intersection point is found by
\begin{equation*}
t=\frac{N_{1}\otimes N_{2}\otimes N_{3} \otimes N_{4}}{\|N_{1}\otimes N_{2}\otimes N_{3} \otimes N_{4}\|}=\left(-\frac{2}{\sqrt{91}},0,-5\sqrt{\frac{2}{91}},-\frac{6}{\sqrt{91}},-\frac{1}{\sqrt{91}}\right)
\end{equation*}
The non-vanishing second fundamental coefficients are\newline
\begin{eqnarray*}&&h_{11}^{1}=-\frac{1}{\sqrt{6}},\hspace{.3cm}h_{22}^{1}=-\frac{1}{\sqrt{6}},\hspace{.3cm}h_{12}^{1}=-\frac{1}{\sqrt{6}},
\hspace{.3cm}h_{11}^{2}=\frac{1}{2\sqrt{2}},\hspace{.3cm}h_{22}^{2}=\frac{1}{\sqrt{2}},\hspace{.3cm}h_{11}^{3}=\frac{2}{3}
,\hspace{.3cm}h_{22}^{3}=\frac{2}{3}\\
&&\hspace{.3cm}h_{12}^{3}=-\frac{2}{3},\hspace{.3cm}h_{13}^{3}=\frac{1}{3}
,\hspace{.3cm}h_{11}^{4}=-\frac{1}{2}\end{eqnarray*}
From the linear system of equations in $(\ref{p18})$, we obtain
\begin{eqnarray*}(u_{1}^{1})^{'}=-0.628971,\hspace{.2cm}(u_{2}^{1})^{'}=0.838628,\hspace{.2cm}(u_{3}^{1})^{'}=-0.209657,\hspace{.2cm}(u_{4}^{1})^{'}=-0.838628\\
(u_{1}^{2})^{'}=0.209657,\hspace{.2cm}(u_{2}^{2})^{'}=-0.419314,\hspace{.2cm}(u_{3}^{2})^{'}=-0.628971,\hspace{.2cm}(u_{4}^{2})^{'}=0.628971\\
(u_{1}^{3})^{'}=0.209657,\hspace{.2cm}(u_{2}^{3})^{'}=-0.419314,\hspace{.2cm}(u_{3}^{3})^{'}=-0.628971,\hspace{.2cm}(u_{4}^{3})^{'}=-1.25794\\
(u_{1}^{4})^{'}=0.419314,\hspace{.2cm}(u_{2}^{4})^{'}=-0.741249,\hspace{.2cm}(u_{3}^{4})^{'}=-0.628971,\hspace{.2cm}(u_{4}^{4})^{'}=-0.209651
\end{eqnarray*}
Thus, we obtain $\kappa_{n}^{1}=-0.879304$, $\kappa_{n}^{2}=0.139867$, $\kappa_{n}^{3}=0.117216$, $\kappa_{n}^{4}=-0.0879121$. Hence, we have
\begin{equation*}\alpha^{''}=-1.321618N_{1}-0.469981N_{2}+0.698466N_{3}-0.285472N_{4}\end{equation*}
Or,\begin{equation*}\alpha^{''}=(-0.839029,-0.658857,0.234991,-0.0995048,0.613453)\end{equation*}
Thus, $\kappa_{1}=\|\alpha^{''}\|=1.25679$\newline
Also the unit normal vector is \newline
n$=(-0.524421,-0.411807, 0.146877, -0.621938, 0.383428)$\newline
From $(\ref{p24})$, we get\newline
$\triangle_{1}=(-1.05495,-1.14286,-0.279735,-0.395605,-0.021978),
\newline\triangle_{2}=(-0.32967,-0.417583,0.248653,0,-0.549451)$\newline
$\triangle_{3}=(-0.32967,-0.417583,-0.404061,0,-0.197802),\;\triangle_{4}=(0,-0.0879121,0,0,0)$\newline
From the linear system of equations in $(\ref{p25})$, we obtain
\begin{eqnarray*}(u_{1}^{1})^{''}=0.106614,\hspace{.2cm}(u_{2}^{1})^{''}=0.161474,\hspace{.2cm}(u_{3}^{1})^{''}=0.89026,\hspace{.2cm}(u_{4}^{1})^{''}=1.05091\\
(u_{1}^{2})^{''}=0.268086,\hspace{.2cm}(u_{2}^{2})^{''}=0.366087,\hspace{.2cm}(u_{3}^{2})^{''}=-0.309769,\hspace{.2cm}(u_{4}^{2})^{''}=-2.69189\\
(u_{1}^{3})^{''}=0.454481,\hspace{.2cm}(u_{2}^{3})^{''}=0.795236,\hspace{.2cm}(u_{3}^{3})^{''}=-0.141793,\hspace{.2cm}(u_{4}^{3})^{''}=0.450137\\
(u_{1}^{4})^{''}=1.67806,\hspace{.2cm}(u_{2}^{4})^{''}=0.2356,\hspace{.2cm}(u_{3}^{4})^{''}=-0.234991,\hspace{.2cm}(u_{4}^{4})^{''}=1.22691
\end{eqnarray*}
Which yields \newline
$\mu_{1}=0.4544,\hspace{.2cm} \mu_{2}=-0.5105,\hspace{.2cm} \mu_{3}=-1.4338, \hspace{.2cm}\mu_{4}=-1.0554$\newline
Then, we have
\begin{equation*}\alpha^{'''}=-(1.25679)^{2}t+1.84188N_{1}+0.451035N_{2}-2.14772N_{3}-1.64788N_{4}\end{equation*}
\begin{equation*}\alpha^{'''}=(2.35545, -1.0554, 0.945301, 0.596496, 0.0935034)\end{equation*}
$\kappa_{1}^{'}=\langle\alpha^{'''}, n\rangle=-0.996915$\newline
Using $\alpha^{'''}$, we obtain
\begin{eqnarray*}(u_{1}^{1})^{'''}=-5.90731,\hspace{.2cm}(u_{2}^{1})^{'''}=2.14186,\hspace{.2cm}(u_{3}^{1})^{'''}=1.51393,\hspace{.2cm}(u_{4}^{1})^{'''}=2.11042\\
(u_{1}^{2})^{'''}=-3.76546,\hspace{.2cm}(u_{2}^{2})^{'''}=2.58181,\hspace{.2cm}(u_{3}^{2})^{'''}=0.76716,\hspace{.2cm}(u_{4}^{2})^{'''}=-2.76882\\
(u_{1}^{3})^{'''}=-0.88196,\hspace{.2cm}(u_{2}^{3})^{'''}=-3.60164,\hspace{.2cm}(u_{3}^{3})^{'''}=-1.14366,\hspace{.2cm}(u_{4}^{3})^{'''}=2.21882\\
(u_{1}^{4})^{'''}=-3.33833,\hspace{.2cm}(u_{2}^{4})^{'''}=0.545801,\hspace{.2cm}(u_{3}^{4})^{'''}=0.596496,\hspace{.2cm}(u_{4}^{4})^{'''}=0.187007
\end{eqnarray*}
Hence, we get
\begin{equation*}\xi_{1}=-21.7826,\hspace{.3cm}\xi_{2}=-3.3273,\hspace{.3cm}\xi_{3}=8.8342,\hspace{.3cm}\xi_{4}=-1.3710\end{equation*}
From $(\ref{p26})$, we have
\begin{eqnarray*}\alpha^{(4)}&=&-3(1.25679)(-0.996915)\left(-\frac{2}{\sqrt{91}},0,-5\sqrt{\frac{2}{91}},-\frac{6}{\sqrt{91}},-\frac{1}{\sqrt{91}}\right)
-41.7712N_{1}\\
&&-29.6654N_{2}+34.1873N_{3}+5.19372N_{4}\end{eqnarray*}
Or,
\begin{equation*}
\alpha^{(4)}=(-30.1443, -1.371, 12.0465, -11.945, 10.9205)\end{equation*}
Thus,
\begin{equation*}b_{3}=\frac{\alpha^{'}\otimes\alpha^{''}\otimes\alpha^{'''}\otimes\alpha^{(4)}}{\|\alpha^{'}\otimes\alpha^{''}\otimes\alpha^{'''}\otimes\alpha^{(4)}\|}
=(0.235522, 0.439366, -0.161734, -0.0297594, 0.851143),\end{equation*}
\begin{equation*}b_{2}=\frac{b_{3}\otimes\alpha^{'}\otimes\alpha^{''}\otimes\alpha^{'''}}{\|b_{3}\otimes\alpha^{'}\otimes\alpha^{''}\otimes\alpha^{'''}\|}
=(0.0469859, -0.0544597, 0.64055, -0.763297, -0.0435799)\end{equation*}
Also,
\begin{equation*}b_{1}=\frac{ b_{2} \otimes b_{3}\otimes \alpha^{'} \otimes \alpha^{''}}{\| b_{2} \otimes b_{3}\otimes \alpha^{'} \otimes \alpha^{''}\|}=(0.648319, -0.724303, -0.161293, -0.0530185, 0.16199)\end{equation*}
Now curvature,
\begin{equation*}\kappa_{2}=\frac{\langle \alpha^{'''},b_{1}\rangle}{\kappa_{1}}=1.68888,\hspace{.5cm}\kappa_{3}=\frac{\langle \alpha^{(4)},b_{2}\rangle}{\kappa_{1}\kappa_{2}}=7.07463\end{equation*}
Also from $(\ref{p4})$ we have, $\kappa_{1}^{''}=35.3284$\newline
Using $\alpha^{(4)}$, we obtain
\begin{eqnarray*}(u_{1}^{1})^{(4)}=44.1618,\hspace{.2cm}(u_{2}^{1})^{(4)}=-21.5916,\hspace{.2cm}(u_{3}^{1})^{(4)}=-4.64151,\hspace{.2cm}(u_{4}^{1})^{(4)}=-16.5865\\
(u_{1}^{2})^{(4)}=22.5702,\hspace{.2cm}(u_{2}^{2})^{(4)}=33.779,\hspace{.2cm}(u_{3}^{2})^{(4)}=-7.12615,\hspace{.2cm}(u_{4}^{2})^{(4)}=55.62\\
(u_{1}^{3})^{(4)}=26.8139,\hspace{.2cm}(u_{2}^{3})^{(4)}=22.2084,\hspace{.2cm}(u_{3}^{3})^{(4)}=-4.81273,\hspace{.2cm}(u_{4}^{3})^{(4)}=-9.77667\\
(u_{1}^{4})^{(4)}=-0.6116,\hspace{.2cm}(u_{2}^{4})^{(4)}=16.2234,\hspace{.2cm}(u_{3}^{4})^{(4)}=-11.945,\hspace{.2cm}(u_{4}^{4})^{(4)}=21.841
\end{eqnarray*}
Hence, we have
\begin{equation*}\eta_{1}=194.9662,\hspace{.3cm}\eta_{2}=-27.5036,\hspace{.3cm}\eta_{3}=-56.7000,\hspace{.3cm}\eta_{4}=29.6466\end{equation*}
Thus we have
\begin{eqnarray*}\alpha^{(5)}&=& \left(-3 (1.25679)^{2} -
    4(1.25679)(35.3284)+(1.25679)^{4}+(1.25679)^{2}(1.68888)^{2}\right)t\\
    &&+313.53N_{1}+145.818N_{2}-210.771N_{3}-46.7969N_{4}\end{eqnarray*}
Or,\begin{equation*}\alpha^{(5)}=(253.719, 29.6467, 57.0616, 143.136, -97.1016)\end{equation*}
Thus, the fourth curvature is given by
\begin{equation*}\kappa_{4}=\frac{\langle \alpha^{(5)},b_{3}\rangle}{\kappa_{1}\kappa_{2}\kappa_{3}}=-1.55521\end{equation*}
Now, to find $\kappa^{1}_{1g}$ and $\tau^{1}_{jg}$, $j=1,2,3$ for hypersurface $M_{1}$,
we have from $(30)$ and $(31)$,
$\kappa^{1}_{1g}=0.584888 ,\tau^{1}_{1g}=-0.774977,\tau^{1}_{2g}=-0.0496875,\tau^{1}_{3g}=-0.276372$.
\noindent Similarly we can find $\kappa^{i}_{1g}$ and $\tau^{i}_{jg}$ for $M_{2},M_{3},M_{4}$.

\end{document}